%% file: Ceretani_et_al._revised_version.tex
\documentclass[graybox]{svmult}

\usepackage{type1cm}        
\usepackage{makeidx}         
\usepackage{graphicx}        
\usepackage{multicol}        
\usepackage[bottom]{footmisc}

\usepackage{newtxtext}       %
\usepackage[varvw]{newtxmath}       

\makeindex             

\begin{document}

\title*{Exact solutions for the fractional 
	nonlinear Boussinesq equation}
\author{A. Ceretani, F. Falcini and R. Garra}
\institute{A. Ceretani \at Department of Mathematics of the Faculty of Exact and Natural Sciences, University of Buenos
	Aires, and Mathematics Research Institute “Luis A. Santal\'o” (IMAS), CONICET, Argentina, \email{aceretani@dm.uba.ar}
\and F. Falcini \at Institute of Marine Sciences, National Research Council (CNR), Rome, Italy, \email{federico.falcini@cnr.it}
\and R. Garra  \at Institute of Marine Sciences, National Research Council (CNR), Rome, Italy, \email{rolinipame@yahoo.it} }
%
%
\maketitle

\abstract{We investigate the existence of exact solutions in closed form to a fractional version of the nonlinear Boussinesq equation for groundwater flow through an unconfined aquifer. We show this fractional equation appears naturally when the classical nonlinear Darcy's law is replaced by a space-fractional one. After a physical discussion on the fractional model, we give several exact solutions in closed form for special choices of initial and boundary data. We provide solutions for steady and unsteady cases, by considering both classical and fractional derivatives in time.}

 \section{Introduction}
This paper is devoted to show that, by simple methods, it is possible to find exact solutions in closed form to a nonlocal generalization of the well-known Boussinesq equation for a horizontal unconfined aquifer,
\begin{equation}\label{Boussinesq}
	\frac{\partial h}{\partial t}(x,t) = k \frac{\partial }{\partial x}\left(h(x,t)\frac{\partial h}{\partial x}(x,t)\right)
+\frac{\Lambda(h(x,t),t)}{n}\qquad x>0,\,t>0,
\end{equation}
where $k,n$ are given positive parameters and $\Lambda$ is a known function. Here, we consider the fully noninear fractional equation
\begin{align}
	\label{Eq-1}\frac{\partial^\gamma h}{\partial t^\gamma}(x,t)=\,&k\frac{\partial}{\partial x}\left(h(x,t)\frac{\partial^\nu h}{\partial x^\nu}(x,t)\right)
	+\frac{\Lambda(h(x,t),t)}{n}\qquad x>0,\,t>0,
\end{align}
with a sink therm $\Lambda$ of the form
\begin{equation}\label{L}
\Lambda(h,t)=-n\phi h,
\end{equation}
where $\phi$ is a non-negative parameter. The operator $\partial^\nu/\partial x^\nu$ is the partial fractional derivative with respect to the space variable $x$, of order $\nu\in(0,1)$, associated to the left-sided Caputo fractional derivative with zero-starting point. The operator $\partial^\gamma/\partial t^\gamma$ is defined analogously if $\gamma\in(0,1)$, whereas it is given by the usual derivative if $\gamma=1$. We recall that the Caputo fractional derivative of order $\alpha\in(0,1)$ is defined by
\begin{equation}\label{caputo}
	\frac{d^\alpha f}{d z^\alpha}(z)=\frac{1}{\Gamma(1-\alpha)}\int_0^z(z-z')^{-\alpha}\frac{d f}{d z'}(z')\,dz'\qquad z\geq 0,
\end{equation}
for absolutely continuous functions $f$ in $[0,\infty)$, where $\Gamma(\cdot)$ is the Gamma-function. Notice that we have a space-time fractional problem if $\gamma\in(0,1)$ and a space-fractional one if $\gamma=1$. In addition, observe that the local 
version ($\gamma=\nu=1$) of equation \eqref{Eq-1} is the usual Boussinesq equation for a horizontal unconfined aquifer.

Our work is motivated by modelling groundwater flows through anisotropic or non-homogeneous soils in hillslope areas. It is widely accepted that the Boussinesq equation is a suitable model to describe groundwater flows through isotropic and homogeneous soils; see for example \cite{brut,DaPo2004,rup} and the references therein. However, in the absence of isotropy or homogeneity and, in particular, under  long-range effects in groudwater flows, the Boussinesq equation can no longer capture the nonlocal effects of presure gradients; see for example \cite{sen,Su,Su1}. In the last decades fractional calculus tools have shown remarkable potential to accurate modeling phenomena in anisotropic, non-homogeneous media, see for example \cite{noi,von}. We show in next Section \ref{Physica} that the space-fractional equation \eqref{Eq-1} arises naturally when the local, nonlinear Darcy's law is replaced by a fractional one in the usual derivation of the Boussinesq equation.

As a first step, we are interested in deriving exact solutions in closed form to our model. This is a hard task because we have to deal with a nonlinear fractional equation. Exact solutions to the Boussinesq equation were given in several works for special choices of initial and boundary data; see for example \cite{DaPo2004,PaEtAl2000,sen}. Following this approach, we consider here a {\em toy model} for which solutions can be obained in closed form. More precisely, we consider the initial and boundary conditions
\begin{align}
	\label{IC-BC-1}h(x,0)=x^\sigma\qquad x\geq 0,\qquad\qquad
	h(0,t)=0\quad t\geq 0,
\end{align}
where $\sigma>0$. 

Frequently, dealing with more realistic fractional models requires to develop methods to find approximate solutions. We refer to the recent monograph \cite{Ka2019} and the references therein for a comprehensive survey on this topic. The exact solutions provided here may serve to test analytical or numerical methods to approximate solutions to fully nonlinear fractional equations. The complexity of fractional models have further motivated the study on generalized solutions to them. We refer, for example, to the recent article \cite{ToRy2020} where the authors investigate the existence and properties of viscosity solutions for a fractional problem obtained from a fractional linear Darcy-type law. To the best of our knowledge, the problem considered here (i.e. the one derived from a fractional {\em nonlinear} Darcy-type law) was not yet investigated in the literature.

The outline of the article is as follows. In Section \ref{Physica} we provide a simple derivation of equation \eqref{Eq-1}, based on the continuity equation and a fractional Darcy-type law. Exact solutions for the steady and unsteady space-fractional problems are given in Sections \ref{Steady} and \ref{Unsteady}, respectively. The time-fractional case is analyzed in Section \ref{TimeFractional}, where exact solutions are provided also for problems with time-dependent coefficients. The paper closes with a breaf discussion in Section \ref{Final}.

\section{Physical motivation}\label{Physica}
The goal of this section is to show that equation \eqref{Eq-1} is a special case of a fractional version of the classical nonlinear Boussinesq equation, which plays a central role to understand groundwater flows through unconfined aquifers.

For isotropic and homogeneous soils, the depth of water $h$ is governed by the continuity equation 
\begin{equation}\label{1}
	n\frac{\partial h}{\partial t}(x,t)= -\frac{\partial q}{\partial x}(x,t)+\Lambda(h(x,t),t),
\end{equation}
where the flow rate $q$ is assumed to be given by the Darcy's law
\begin{equation}\label{2}
	q(x,t)= -h(x,t) K_s \sin (\theta) -h(x,t) \frac{\partial h}{\partial x}(x,t)K_s \cos(\theta),
\end{equation}
see for example \cite{bartlett}. Here, $x$ is the coordinate parallel to the bed slope and $\theta$ is the angle of the bed slope to the horizontal. The parameter $n$ is the soil porosity, $K_s$ is the saturated hydraulic conductivity, and $\Lambda(h,t)$ is a source or sink term. 

Replacing equation \eqref{2} into equation \eqref{1} we directly obtain the nonlinear Boussinesq equation
\begin{equation}\label{3}
	\frac{\partial h}{\partial t}(x,t) = k \frac{\partial }{\partial x}\left(h(x,t)\frac{\partial h}{\partial x}(x,t)\right)+v_s\frac{\partial h}{\partial x}(x,t)+\frac{\Lambda(h(x,t),t)}{n},
\end{equation}
where we have written $v_s = (K_s/n) \sin (\theta)$ and $k = (K_s/n) \cos (\theta)$ for shortness.

In the absence of homogeneity or isotropy in the porous medium, a model based on the local Darcy's law \eqref{2} may no longer capture the effects of nonlocal pressure differences. The role of them on the Darcy's law was firstly investigated by Sen and Ramos in \cite{sen}, for a steady problem. There, the porous medium is conceived as a network of short, medium and long-range interstitial channels with impermeable walls, where the flow is driven by both local and nonlocal pressure gradients. The flow rate $q$ is then given by 
\begin{equation}\label{general}
q(x) = -\int_{-\infty}^{+\infty} f(x,x')(p(x')- p(x)) dx',
\end{equation}
where $p$ is the pressure field and $f$ is a nonlocal flow conductivity. The later takes into account the effects of the channel lenght and thickness between the locations $x'$ and $x$, and it is directly related to the medium properties. Both the classical Darcy's law \eqref{2} and some nonlocal versions of it can be obtained as special cases of \eqref{general}; see \cite{sen} for the details. In particular, by considering a power-law nonlocal conductivity $f$ one finds a fractional Darcy-type law given in terms of the Riemann-Liouville fractional derivative operator \cite{sen}. Encouraged by this pioneering work, models based on fractional Darcy-type laws were investigated in several papers, see for example \cite{physica} and the references therein. In the same vein of them, we consider here that the flow rate $q$ is given by the nonlinear fractional Darcy-type law
\begin{equation}\label{2bis}
q(x,t) = -h(x,t) K_s \sin (\theta) -  h(x,t) \frac{\partial^\nu h}{\partial x^\nu}(x,t)K_s \cos(\theta), 
\end{equation}  
where $\nu\in(0,1)$. Replacing \eqref{2bis} into the continuity equation \eqref{1} we obtain the nonlinear space-fractional Boussinesq equation
\begin{equation}\label{4}
	\frac{\partial h}{\partial t}(x,t)= k \frac{\partial}{\partial x}\left(h \frac{\partial^\nu h}{\partial x^\nu}(x,t)\right)+v_s\frac{\partial h}{\partial x}(x,t)+\Lambda(h(x,t),t).
\end{equation}

To keep the problem as simple as possible while retaining the nonlinearity, we consider a horizontal aquifer (then $\theta=0$ and so $v_s$) and a sink term of the form $\Lambda(h,t) = -n\phi h$, where $\phi$ is a transfer coefficient that controls the magnitude of the seepage loss of water into the bedrock; see \cite{bartlett}. Thus, equation \eqref{4} reduces to \eqref{Eq-1} with $\gamma=1$.

\section{The steady solution}\label{Steady}
We first consider the steady case, which is given by
the boundary value problem
\begin{equation}
	k\frac{d}{dx}\left(h(x)\frac{d^\nu h}{dx^\nu}(x)\right) = \phi h(x) \qquad x>0,\qquad\qquad h(0)=0,
\end{equation}
for a nonlinear fractional ordinary differential equation. This problem admits a solution in the power-law form 
$h(x) = cx\,^\beta$ provided that $c$ and $\beta$ are properly chosen, as we show below. Notice that any function $h$ in this form satisfies the boundary condition $h(0)=0$. 

Exploiting that
\begin{equation}\label{dr}
	\frac{d^\nu}{d x^\nu}(x\,^\beta)=\frac{\Gamma(\beta+1)}{\Gamma(\beta+1-\nu)}x\,^{\beta-\nu},
\end{equation}
we obtain that $h(x) = cx\,^\beta$ satisfies the equation if and only if
\begin{equation}\label{cond}
k c(2\beta-\nu)\frac{\Gamma(\beta+1)}{\Gamma(\beta+1-\nu)}x^{2\beta-\nu-1}={\phi x^\beta}.
\end{equation}

If $\phi=0$, equation \eqref{cond} is fulfilled only if $c=0$. If $\phi>0$, an elementary computation shows that \eqref{cond} holds true only if $\beta=\nu+1$ and $c=\phi/k\Gamma(\nu+3)$. Thus, we find the exact solution
\begin{equation}
h(x) =\frac{{\phi x^{\nu+1}}}{k\Gamma(\nu+3)} \qquad x\geq 0.
\end{equation}

\begin{remark}
Following the same argument, we find the exact solution 
\begin{equation}
	h_\ell(x)=\frac{{\phi x^2}}{6k} \qquad x\geq 0,
\end{equation}
for the associated nonlinear local problem
\begin{equation}
	k\frac{d}{dx}\left(h(x)\frac{d h}{dx}(x)\right) = {\phi h(x)}\qquad x>0,\qquad\qquad h(0)=0.
\end{equation}
\end{remark}

\section{The unsteady space-fractional case}\label{Unsteady}
We now address the unsteady space-fractional problem given by
\begin{equation}\label{space}
\frac{\partial h}{\partial t}(x,t) = k\frac{\partial}{\partial x}\left(h(x,t)\frac{\partial^\nu h}{\partial x^\nu}(x,t)\right) - \phi h(x,t)\qquad x>0,\,t>0,
\end{equation}
subject to the initial and boundary conditions \eqref{IC-BC-1}.

Motivated by the steady solution found in Section \ref{Steady}, we assume that the exponent $\sigma$ in the initial data $h(x,0)=x^{\sigma}$ is given by $\sigma=\nu+1$. We show below that this problem admits a similarity-type solution in the form
\begin{equation}\label{h}
h(x,t) = x^{\nu+1}f(t),
\end{equation}
for a suitable selection of the function $f$. 

Exploiting formula \eqref{dr} again, we observe that $h(x,t) = x^{\nu+1}f(t)$ solves equation \eqref{space} if and only if $f$ solves the following boundary value problem for a nonlinear ordinary differential equation: 
\begin{equation}\label{f}
	\frac{df}{dt}(t) = k\,\Gamma(\nu+3) f(t)^2-\phi f(t)\qquad t>0,\qquad f(0)=1.
\end{equation}

If $\phi=0$, a simple computation shows that 
\begin{equation*}
f(t)=\frac{1}{1-k\Gamma(\nu+3)t},
\end{equation*}
solves \eqref{f} for $t\geq 0$ with $t\neq 1/k\Gamma(\nu+3)$. If $\phi>0$, we find that 
\begin{equation*}
f(t)=\frac{\phi}{e^{\phi t}(\phi - k \Gamma(\nu+3 ))+ k \Gamma(\nu+3)},
\end{equation*}
solves \eqref{f} for $t\geq 0$ if $\phi\geq k\Gamma(\nu+3)$, or for $t\geq 0$ with $t\neq\frac{1}{\phi}\ln\left(\frac{k \Gamma(\nu +3 )}{k \Gamma(\nu +3)-\phi}\right)$ if $0<\phi<k\Gamma(\nu+3)$. Phisically realistic solutions for $\phi=0$ must be non-negative. This is achieved only if $0\leq t\leq \frac{1}{k\Gamma(\nu+3)}$. If $\phi>0$, physically realistic solutions must be, in addition, non-increasing with respect to time at every $x>0$. This occurs only if $\phi\geq k\Gamma(\nu+3)$. Thus, we obtain the exact solution
\begin{equation*}
h(x,t)=\left\{\begin{array}{cc}
\displaystyle\frac{x^{\nu+1}}{1-k\,\Gamma(\nu+3)t}&\quad x\geq 0,\,0\leq t<\frac{1}{k\Gamma(\nu+3)},\,\,\,\mathrm{if}\quad\phi=0,\\[0.5cm]
\displaystyle\frac{\phi \ x^{\nu+1}}{e^{\phi t}(\phi - k \Gamma(\nu+3 ))+ k \Gamma(\nu+3)}&\quad x\geq 0,\, t\geq 0,\,\,\,\mathrm{if}\quad\phi\geq k\Gamma(\nu+3).
\end{array}\right.
\end{equation*}
%

\begin{remark}
As for the steady case, similar steps allow us to obtain the exact solution 
\begin{equation*}
h_{\ell}(x,t)=\left\{\begin{array}{clcc}
\displaystyle \frac{x^2}{1-6kt}&\qquad x\geq 0,\,0\leq t<\frac{1}{6k}&\quad\mathrm{if}\quad&\phi=0,\\[0.5cm]
\displaystyle \frac{\phi \ x^{2}}{e^{\phi t}(\phi - 6k)+ 6k}&\qquad x\geq 0,\, t\geq 0&\quad\mathrm{if}\quad&\phi\geq 6k,
\end{array}\right.
\end{equation*}
to the associated local problem with initial data $h(x,0)=x^{2}$.
\end{remark}

Minor changes on the above arguments allow to find exact solutions in the form
\begin{equation}\label{hbis}
h(x,t)=x^{\nu+1}f_1(t)+f_2(t),
\end{equation}
where $f_1$ and $f_2$ satisfy the system of nonlinear ordinary differential equations
\begin{align*}
\frac{df_1}{dt}(t)=\,&k\Gamma(\nu+3)f_1(t)^2-\phi f_1(t),\\
\frac{dg}{dt}(t)=\,&f_1(t)f_2(t)\Gamma(\nu+2)-\phi f_2(t),
\end{align*} 
with the initial conditions $f_1(0)=1$, $f_2(0)=0$.

The strategy developed so far to obtain exact solutions in the form \eqref{h} or \eqref{hbis} can be framed into the invariant subspace method (see for example \cite{gala} for the general theory and \cite{gazi} for applications to nonlinear fractional differential equations). We briefly recall the main idea of this method: consider the equation
\begin{equation}\label{pro}
\frac{\partial h}{\partial t}(x,t)= F[h](x,t),
\end{equation}
where $F[h]$ is a nonlinear operator. A space of functions $W$ is said to be {\em invariant} under the operator $F$ if $F[h]\in W$ for every $h\in W$. It is straightforward that if $W$ admits a finite basis $\{w_1,\hdots,w_n\}$, then $W$ is invariant under $F$ if and only if
there exist $n$ scalar functions $\Phi_1,\hdots,\Phi_n$ defined on $\mathbb{R}^n$ that satisfy
\begin{equation*}
F[\sum_{k=1}^nf_kw_k]=\sum_{k=1}^n\Phi_k(f_1,\hdots,f_n)w_k,
\end{equation*}
for every $(f_1,f_2,\hdots,f_n)\in\mathbb{R}^n$. Thus, if $W$ is as before, then equation \eqref{pro} admits a solution in the form
\begin{equation*}
h(x,t)=\sum_{k=1}^nf_k(t)w_k(x),
\end{equation*}  
if and only if $f_1,\hdots,f_n$ satisfy 
\begin{align*}
\frac{df_1}{dt}(t)=\,&\Phi_1(f_1(t),\hdots,f_n(t)),\\[0.25cm]
\frac{df_2}{dt}(t)=\,&\Phi_2(f_1(t),\hdots,f_n(t)),\\
\vdots\\
\frac{df_n}{dt}(t)=\,&\Phi_n(f_1(t),\hdots,f_n(t)).
\end{align*} 
In this way, we reduce the original problem \eqref{pro} to a system of ordinary differential equations that can be solved exactly in many cases. The useful role of this method to find exact solutions for nonlinear fractional partial differential equations has been shown in recent publications, we refer for example to \cite{saha,prak}. 

Lets now go back to our problem. The operator $F$ associated to the equation \eqref{space} is given by
\begin{equation}\label{F}
F[h]=k\frac{\partial }{\partial x}\left(h\frac{\partial^\nu h}{\partial x^\nu}\right)-\phi h.
\end{equation}
Let $w_1$ and $w_2$ be the functions given by
\begin{equation}\label{w12}
w_1(x)=x^{\nu+1}\qquad\text{and}\qquad w_2(x)=1.
\end{equation}
A straighforward computation shows that the spaces $W_1$ and $W_2$ spanned by $\{w_1\}$ and $\{w_1,w_2\}$, respectively, are invariant under the operator $F$. The solutions in the form \eqref{h} are those related to $W_1$, and the solutions in the form \eqref{hbis} are those related to $W_2$.

Exploiting the invariant subspace method, one can look for exact solutions to equation \eqref{space} by first finding a finite dimensional invariant subspace for the operator $F$ given in \eqref{F}. To illustrate, it is also possible to find exact solutions from the invariant subspaces $\tilde{W}_1$ and $\tilde{W}_2$ spanned by $\{\tilde{w}_1\}$ and $\{\tilde{w}_1,w_2\}$, respectively, where 
\begin{equation}\label{w1t}
\tilde{w}_1(x)=x^{\nu/2},
\end{equation}
provided that the exponent in the initial data $h(x,0)=x^{\sigma}$ in \eqref{IC-BC-1} is $\sigma=\nu/2$. For example, the solution related to $\tilde{W}_1$ is given by   
\begin{equation}
h(x,t) =  x^{\nu/2} e^{-\phi t}\qquad x\geq 0,\,t\geq 0.
\end{equation}
We omit the details to avoid repetition.

\section{The time-fractional case}\label{TimeFractional}

We finally consider the space-time fractional problem given by equation \eqref{Eq-1} with $\gamma\in(0,1)$ and the sink term given by \eqref{L}, with the initial and boundary conditions \eqref{IC-BC-1}.

The invariant subspace method described at the end of Section \ref{Unsteady} extends naturally to equations in the form
\begin{equation}
\frac{\partial^\gamma h}{\partial t^\gamma}=F[h].
\end{equation}

Let $W_1$, $W_2$, $\tilde{W}_1$, $\tilde{W}_2$ be the subspaces spanned by $\{w_1\}$, $\{w_1,w_2\}$, $\{\tilde{w}_1\}$, $\{\tilde{w}_1,w_2\}$, where $w_1$, $w_2$ and $\tilde{w}_1$ are given in \eqref{w12} and \eqref{w1t} (see Section \ref{Unsteady}). We already know that they are invariant spaces under the operator $F$ given by \eqref{F}. This motivates to look for solutions that belong to some of them. For example, by considering the invariant subspace $W_1$, we find the exact solution 
\begin{equation}
h(x,t) =  \ x^{\nu/2} E_\gamma(-\phi\, t^\gamma)\qquad x\geq 0,\,t\geq 0,
\end{equation}
for the case when the exponent in the initial data is $\sigma=\nu/2$. Here, 
\begin{equation}
	E_\gamma (-\phi\, t^\gamma) = \sum_{i=0}^\infty\frac{(-\phi\, t^\gamma)^i}{\Gamma(\gamma i +1)},
\end{equation}
is the one-parameter Mittag-Leffler function (see \cite{main}).

\begin{remark}The same arguments enable to find exact solutions when the transfer coefficient $\phi$ depends on time according to $\phi(t)=\lambda\,t^{\beta}$, where $\lambda >0$ and $\beta\in\mathbb{R}$.
Here we recall that the so-called Kilbas-Saigo function (see \cite{main} and \cite{bo} for details about this special function)
\begin{align}\label{saigo}
	E_{a,1+\frac{b}{a},\frac{b}{a}} & \left(-
	\lambda t^{a + b} \right)
	= 1 + \sum_{i=1}^\infty (-\lambda)^i t^{i \left( a + b \right)}
	\prod_{j=0}^{i-1} \frac{\Gamma \left( a \left( j+j \frac{b}{a}
		+ \frac{b}{a} \right)+1 \right)}{\Gamma \left(
		a \left( j + j \frac{b}{a} + \frac{b}{a} 
		+ 1 \right) +1 \right)},       
\end{align}

solves the fractional Cauchy problem:
\begin{equation} \label{odefr}
	\begin{cases}
		\frac{d^{a}y}{dt} \left(t \right) = -\lambda t^b y \left( t \right), &
		t \ge 0, \: a \in \left( 0, 1 \right], \: -a < b 
		\leq 1-a, \\
		y \left( 0 \right) = 1.
	\end{cases}
\end{equation}

Therefore, in this case we obtain the solution
\begin{equation}
	h(x,t) =  x^{\nu/2} E_{\gamma,1+\frac{\beta}{\gamma},\frac{\beta}{\gamma}} \left(-
	\lambda\, t^{\beta + \gamma} \right)\qquad x\geq 0,\,t\geq 0,
\end{equation}
under the assumptions $\beta\in(-\gamma,1-\gamma]$ and $\sigma=\nu/2$.
\end{remark}

\section{Conclusions}\label{Final} 

We have explored the existence of exact solutions in closed form to a fully nonlinear fractional equation involving Caputo partial derivatives. This equation is closely related to the Boussinesq equation for unconfined horizontal aquifers and it was obtained by a fractional Darcy-type law for the flow rate. Supplemented with appropriated initial and boundary conditions, the fractional equation considered here may serve as a model for groundwater flows through anisotropic or non-homogeneous porous medium. As a first step to validate the later, we addressed here a toy model that allowed us to find exact solutions in closed form. This model was obtained by a specific selection of initial and boundary conditions. We provided exact solutions to steady and unsteady space-fractional problems, with and without memory effects. These solutions may be used to test analytical or numerical methods to approximate solutions to more realistic, steady and unsteady models, based on the fractional equation proposed here.  

Nonlinearity, nonlocality, and memory effects may be directly related to unknown, long-range factors that vary across landscape types. This represents a core general  hydrologic descriptions at both hillslope and watershed scales, since nonlinearity, nonlocality and memory may affect the 
groundwater flow at different time and spatial scales.
As for \cite{bartlett}, our general solutions also have an explicit dependence on the along-channel variable $x$, a feature that allows for analytic investigations of surface topography and water table profile. Finally, our set of solutions reinforce previous finding on fractional studies on the Boussinesq equation for groundwater flow (\cite{Su,Su1}), presenting new solutions related to different initial and boundary conditions and thus providing new insights into how water flows in fractal media, also for transient problems.


%
\begin{acknowledgement}
This work was partially supported by EMODnet (European Marine Observation and Data Network Physics) and by CNR, in the frame of the Italian Flag Project Ritmare and the Italian Marine Strategy Framework Directive Programmes. \\
This is a preprint of the following chapter: A. Ceretani, F. Falcini and R. Garra, Exact solutions for the fractional nonlinear
Boussinesq equation, published in "Fractional Differential Equations: Modeling, Discretization, and Numerical Solvers", edited by: Cardone A., Donatelli M., Durastante F., Garrappa R., Mazza M., Popolizio M., 2023, Springer INdAM Series
\end{acknowledgement}
%

\input{references}

\end{document}

%% file: references.tex
%
%
%